\documentclass[11pt]{amsart}

\usepackage{amscd, amsfonts, mathrsfs, amsmath, amssymb, amsthm, mathtools}
\usepackage{xcolor}

\newtheorem{theorem}{Theorem}[section]

\theoremstyle{definition}

\newtheorem{remark}[theorem]{Remark}

\def \Q {\mathbb{Q}}
\def \Qbar {\overline{\mathbb{Q}}}
\DeclareMathOperator{\rk}{rk}
\DeclareMathOperator{\Cl}{Cl}
\DeclareMathOperator{\Jac}{Jac}

\DeclareMathOperator{\tors}{tors}

\numberwithin{equation}{section}

\title[Imaginary Quadratic Fields with Class Group of $5$-rank at least $2$]{Counting Imaginary Quadratic Fields with an Ideal Class Group of $5$-rank at least $2$}
\author{Kollin Bartz}
\author{Aaron Levin}
\author{Aman Dhruva Thamminana}

\address{
}
\address{
Department of Mathematics\\
Michigan State University\\
East Lansing, MI 48824\\
USA
}
\email{bartzkol@msu.edu}
\email{levina@msu.edu}
\email{thammina@msu.edu}

\date{}

\begin{document}

\begin{abstract}
We prove that there are $\gg\frac{X^{\frac{1}{3}}}{(\log X)^2}$ imaginary quadratic fields $k$ with discriminant $|d_k|\leq X$ and an ideal class group of $5$-rank at least $2$. This improves a result of Byeon, who proved the lower bound $\gg X^{\frac{1}{4}}$ in the same setting. We use a method of Howe, Lepr\'{e}vost, and Poonen to construct a genus $2$ curve $C$ over $\mathbb{Q}$ such that $C$ has a rational Weierstrass point and the Jacobian of $C$ has a rational torsion subgroup of $5$-rank $2$. We deduce the main result from the existence of the curve $C$ and a quantitative result of Kulkarni and the second author.  
\end{abstract}

\maketitle 

\section{Introduction}

Given an integer $m>1$, it has been known since Nagell's 1922 result \cite{Nagell} that there are infinitely many imaginary quadratic number fields with class number divisible by $m$; the analogous result for real quadratic fields was proved in the early 1970's independently by Yamamoto \cite{Yam} and Weinberger \cite{Wein}. Quantitative results giving a lower bound for the number of such fields were given later by Murty \cite{Murty}, Soundararajan \cite{Sound}, and Yu \cite{Yu}.

More generally, one can study the $m$-rank of the ideal class group. If $A$ is a finitely generated abelian group, we define the $m$-rank of $A$, $\rk_m A$, to be the largest integer $r$ such that $A$ has a subgroup isomorphic to $(\mathbb{Z}/m\mathbb{Z})^r$. For a number field $k$, we let $\Cl(k)$ denote its ideal class group, and let $d_k$ denote its (absolute) discriminant. Concentrating on the imaginary quadratic case, we let $\mathcal{N}^-(m^r;X)$ denote the number of imaginary quadratic number fields $k$ satisfying $|d_k|\leq X$ and  $\rk_m\Cl(k)\geq r$.

With this notation, improving on the results of Murty \cite{Murty}, Soundararajan \cite{Sound} proved (for any $\epsilon>0$),
\begin{align*}
\mathcal{N}^-(m;X)\gg 
\begin{cases}
X^{\frac{1}{2}+\frac{2}{m}-\epsilon} & \quad m\equiv 0\pmod{4},\\
X^{\frac{1}{2}+\frac{3}{m+2}-\epsilon} & \quad m\equiv 2\pmod{4}.
\end{cases}
\end{align*}

Since trivially $\mathcal{N}^-(m;X)\geq \mathcal{N}^-(2m;X)$, one also finds inequalities for odd $m$. For $m=3$, better inequalities have been proved by Heath-Brown \cite{HB} and Lee, Lee, and Yoo \cite{LLY}.

For rank $r=2$, we have the following quantitative result of Kulkarni and the second author, slightly improving an earlier result of Byeon \cite{Byeon}.

\begin{theorem}[Kulkarni-Levin \cite{KL}]
\label{chyp2}
Let $m>1$ be an integer.  Then
\begin{align*}
\mathcal{N}^-(m^2;X)&\gg X^{\frac{1}{m}}/(\log X)^2.
\end{align*}
\end{theorem}

When $m$ is odd, a better lower bound seems to be known only for $m=3,5,7$. Improving on results of Luca and Pacelli \cite{LP} and Yu \cite{Yu}, Lee, Lee and Yoo \cite{LLY} proved
\begin{align*}
\mathcal{N}^{-}(3^2;X)\gg X^{\frac{2}{3}}.
\end{align*}
For $m=5,7$, Byeon proved
\begin{theorem}[Byeon \cite{Byeon3}]
We have
\begin{align*}
\mathcal{N}^-(5^2;X)&\gg X^{\frac{1}{4}}, \\
\mathcal{N}^-(7^2;X)&\gg X^{\frac{1}{4}}.
\end{align*}
\end{theorem}

The goal of this note is the following improvement to Byeon's lower bound for $\mathcal{N}^-(5^2;X)$:
\begin{theorem}
\label{BLT}
We have
\begin{align*}
\mathcal{N}^-(5^2;X)&\gg\frac{X^{\frac{1}{3}}}{(\log X)^2}.
\end{align*}
\end{theorem}

The proof is based on the ``geometric approach" to ideal class groups, which was formalized in work of the second author \cite{LevBordeaux} and Gillibert and the second author \cite{GL1}. More precisely, we use a recent enumerative result of Kulkarni and the second author (building on the work of \cite{GL1}):

\begin{theorem}[Kulkarni-Levin \cite{KL}]
\label{chyp}
Let $C$ be a smooth projective hyperelliptic curve over $\Q$ with a $\Q$-rational Weierstrass point. Let $g=g(C)$ denote the genus of $C$ and let $\Jac(C)(\Q)_{\rm{tors}}$ denote the rational torsion subgroup of the Jacobian of $C$.  Let $m>1$ be an integer and let
\begin{align*}
r=\rk_m \Jac(C)(\Q)_{\rm{tors}}.
\end{align*}
Then
\begin{align*}
\mathcal{N}^-(m^r;X)&\gg \frac{X^{\frac{1}{g+1}}}{(\log X)^2}.
\end{align*}
\end{theorem}

In view of Theorem \ref{chyp}, in order to prove Theorem \ref{BLT} we only need to find a genus $2$ curve $C$ over $\Q$ with a $\Q$-rational Weierstrass point and
\begin{align*}
\rk_5 \Jac(C)(\Q)_{\rm{tors}}\geq 2.
\end{align*}

In fact, in this case we must have equality $\rk_5\Jac(C)(\Q)_{\rm{tors}}=2$, since it is well known from a Weil pairing argument that if $k$ does not contain an $m$th root of unity, then $\rk_m \Jac(C)(k)_{\rm{tors}}\leq g(C)$. To find a suitable genus $2$ curve $C$ we examine a construction of Howe, Lepr\'{e}vost, and Poonen \cite{HLP}, who studied the problem of constructing large torsion subgroups of Jacobians of curves of genus $2$ and $3$. Indeed, one of the aims of the ``geometric approach" to ideal class groups is to take advantage of such constructions for rational torsion in Jacobians, and to leverage these arithmetic-geometric results to study ideal class groups of number fields. 

\section{A construction of Howe, Lepr\'{e}vost, and Poonen}

We will construct a genus $2$ curve $C$ over $\mathbb{Q}$ with a rational Weierstrass point that admits two independent maps of degree $2$ to elliptic curves $E_1$ and $E_2$, with each elliptic curve possessing a rational $5$-torsion point. In this case, the maps $C\to E_1$ and $C\to E_2$ induce an isogeny between $\Jac(C)$ and $E_1\times E_2$ of degree coprime to $5$, and it follows that $\rk_5 \Jac(C)(\Q)_{\rm{tors}}\geq 2$ (in fact, we have equality as noted above).

For this purpose, we use a construction of Howe, Lepr\'{e}vost, and Poonen \cite{HLP} to produce many candidate curves $C$, with the exception that they may not possess a rational Weierstrass point. We then use a brute force search among the produced curves to find a curve $C$ with a rational Weierstrass point.

We begin by summarizing the relevant facts. Recall that elliptic curves with a rational $10$-torsion point form a $1$-parameter family. Indeed, from \cite[Table 6]{HLP}, the universal elliptic curve over $X_1(10)$ can be given explicitly as
\begin{multline}
\label{universal}
E_t:y^2=x(x^2-(2t^2-2t+1)(4t^4-12t^3+6t^2+2t-1)x\\
+16(t^2-3t+1)(t-1)^5t^5),
\end{multline}
and from \cite[Table 7]{HLP} a point of order $10$ on $E_t$ is
\begin{align*}
(x,y)=(4(t-1)(t^2-3t+1)t^3,4(t-1)(t^2-3t+1)t^3(2t-1)).
\end{align*}
The discriminant of $E_t$ modulo squares is
\begin{align*}
\Delta_{10}(t)=(2t-1)(4t^2-2t-1).
\end{align*}

Let $E_t$ and $E_u$ be two elliptic curves over $\mathbb{Q}$ with a rational $10$-torsion point corresponding to the choices $t,u\in\mathbb{Q}$. From \cite[Proposition 3]{HLP}, one can find a genus $2$ curve $C$ with Jacobian $\Jac(C)$ that is $(2,2)$-isogenous to $E_t\times E_u$ if there is an isomorphism of Galois modules $E_t[2](\Qbar)$ to $E_u[2](\Qbar)$ that does not come from an isomorphism of elliptic curves (see \cite[Proposition 3]{HLP} for the precise statement). As noted in \cite{HLP}, since  such an isomorphism must map the rational $2$-torsion point to the rational $2$-torsion point, the isomorphism exists if and only if the quadratic field defined by the non-rational $2$-torsion points on $E_t$ and $E_u$ are equal, and this occurs if and only if the discriminants of $E_t$ and $E_u$ are equal modulo squares. Then the problem of finding suitable elliptic curves $E_t$ and $E_u$ reduces to finding rational solutions to the equation
\begin{align} 
\label{disceqn2} 
\Delta_{10}(t)z^2=\Delta_{10}(u).
\end{align} 
In stating an explicit result, it will be slightly more convenient to use the equivalent equation
\begin{align}
\label{disceqn} 
(2t-1)^5(4t^2-2t-1)z^2=(2u-1)^5(4u^2-2u-1),
\end{align} 
coming from considering the discriminant of the quadratic factor on the right-hand side of \eqref{universal}.

Then given a solution to \eqref{disceqn}, we have the following explicit result constructing a genus $2$ curve $C$ with Jacobian $\Jac(C)$ that is $(2,2)$-isogenous to $E_t\times E_u$.

\begin{theorem}
\label{thmain}
Let $(t,u,z)\in\mathbb{Q}^3$ be a solution to \eqref{disceqn}. Let
\small
\begin{align*}
a=&2(8t^6z - 32t^5z + 40t^4z - 20t^3z + 4tz - 8u^6 + 32u^5 - 40u^4 + 20u^3 - 4u - z + 1),\\
a_0=&-64t^5(t-1)^5(t^2-3t+1),\\
a_2=&2(32t^{12}z - 256t^{11}z + 832t^{10}z - 1440t^9z + 1440t^8z - 960t^7z + 64t^6u^6 - 256t^6u^5 \\
&+    320t^6u^4 - 160t^6u^3 + 32t^6u + 640t^6z - 8t^6 - 256t^5u^6 + 1024t^5u^5 - 1280t^5u^4 \\
&+
    640t^5u^3 - 128t^5u - 480t^5z + 32t^5 + 320t^4u^6 - 1280t^4u^5 + 1600t^4u^4 - 800t^4u^3 \\
		&+
    160t^4u + 240t^4z - 40t^4 - 160t^3u^6 + 640t^3u^5 - 800t^3u^4 + 400t^3u^3 - 80t^3u - 40t^3z\\
		&+
    20t^3 - 16t^2z + 32tu^6 - 128tu^5 + 160tu^4 - 80tu^3 + 16tu + 8tz - 4t - 8u^6 + 32u^5\\
		&-
    40u^4 + 20u^3 - 4u - z + 1),
\\
a_4=&384t^7z^2 - 128t^6u^6z + 512t^6u^5z - 640t^6u^4z + 320t^6u^3z - 64t^6uz - 1152t^6z^2 + 16t^6z\\
& + 
    512t^5u^6z - 2048t^5u^5z + 2560t^5u^4z - 1280t^5u^3z + 256t^5uz + 1344t^5z^2 - 64t^5z\\
		&-
    640t^4u^6z + 2560t^4u^5z - 3200t^4u^4z + 1600t^4u^3z - 320t^4uz - 720t^4z^2 + 80t^4z\\
		&+
    320t^3u^6z - 1280t^3u^5z + 1600t^3u^4z - 800t^3u^3z + 160t^3uz + 120t^3z^2 - 40t^3z + 48t^2z^2\\
    &- 64tu^6z + 256tu^5z - 320tu^4z + 160tu^3z - 32tuz - 24tz^2 + 8tz - 64u^{12} + 512u^{11}\\
		&-
    1664u^{10} + 2880u^9 - 2880u^8 + 1536u^7 + 16u^6z - 128u^6 - 64u^5z - 384u^5 + 80u^4z\\
		&+ 240u^4 -
    40u^3z - 40u^3 - 16u^2 + 8uz + 8u + 3z^2 - 2z - 1,
\\
a_6=&z(-128t^7z^2 + 384t^6z^2 - 448t^5z^2 + 240t^4z^2 - 40t^3z^2 - 16t^2z^2 + 8tz^2 + 64u^{12} \\
&- 512u^{11} +
    1664u^{10} - 2880u^9 + 2880u^8 - 1536u^7 + 128u^6 + 384u^5 - 240u^4 + 40u^3\\
		&+ 16u^2 -
    8u - z^2 + 1).
\end{align*}
\normalsize
Suppose that $a\neq 0$ (or equivalently $z\neq \frac{8u^6 - 32u^5 + 40u^4 - 20u^3 + 4u - 1}{8t^6 - 32t^5 + 40t^4 - 20t^3 + 4t - 1}$) and
\begin{align*}
t,u\notin \left\{0,\frac{1}{2},1\right\}.
\end{align*}
Then
\begin{align*}
y^2=a(a_6x^6+a_4x^4+a_2x^2+a_0)
\end{align*}
defines a genus $2$ curve $C$, $\Jac(C)$ is $(2,2)$-isogenous to $E_t\times E_u$, and in particular, $\rk_5\Jac(C)(\mathbb{Q})_{\tors}=2$. Explicitly, $E_t$ is isomorphic to the elliptic curve
\begin{align*}
E_t':y^2=a(a_0x^3+a_2x^2+a_4x+a_6),
\end{align*}
$E_u$ is isomorphic to the elliptic curve
\begin{align*}
E_u':y^2=a(a_6x^3+a_4x^2+a_2x+a_0),
\end{align*}
and the isogeny may be induced by the two bielliptic maps
\begin{align*}
C&\to E_t'\\
(x,y)&\mapsto (1/x^2,y/x^3)
\end{align*}
and
\begin{align*}
C&\to E_u'\\
(x,y)&\mapsto (x^2,y),
\end{align*}
respectively.
\end{theorem}

Theorem \ref{thmain} is derived from Proposition 4 of \cite{HLP}. However, once the computation is made, it is possible to give an independent direct (computational) proof of the result, as we now describe.

\begin{proof}
We first verify that $y^2=a(a_6x^6+a_4x^4+a_2x^2+a_0)$ defines a genus $2$ curve $C$, or equivalently that $a(a_6x^6+a_4x^4+a_2x^2+a_0)$ is a separable sextic polynomial. By assumption, $a\neq 0$. We have $a_6=zf(t,u,z^2)$ for a certain polynomial $f$. Using \eqref{disceqn} to eliminate $z^2$ in the polynomial $f$ and factoring, we find that $a_6\neq 0$ if $t,u\notin \left\{0,\frac{1}{2},1\right\}
$ and neither $t$ nor $u$ are roots of $4x^2-2x-1$ or $x^2-3x+1$; in fact, this last condition is vacuous since $t,u\in \mathbb{Q}$ and the two quadratic polynomials are irreducible over $\mathbb{Q}$. Therefore $a(a_6x^6+a_4x^4+a_2x^2+a_0)$ is a sextic polynomial. Finally, a direct computation gives that the irreducible factors (in $\mathbb{Q}[t,u,z]$) of the discriminant of $a_6x^6+a_4x^4+a_2x^2+a_0$ are $a, a_6, t, t-1, 2t-1, t^2-3t+1$, and $4t^2-2t-1$. Since none of these factors vanishes under our hypotheses, the sextic polynomial is separable.

It is well-known, going back to Jacobi in the 19th century, that $\Jac(C)$ is $(2,2)$-isogenous to $E_t'\times E_u'$, induced by the given maps (see \cite[Theorem 14.1.1]{CF}). Finally, by computing reduced Weierstrass forms of the curves, one can check by direct computation that $E_t$ and $E_{t'}$ are isomorphic, and that $E_u$ and $E_u'$ are isomorphic (for the latter isomorphism, one must compute modulo the relation \eqref{disceqn}).

\end{proof}

\begin{remark}
\label{2rem}
In fact, it follows from \cite{HLP} that $\Jac(C)$ contains a rational torsion subgroup isomorphic to $\mathbb{Z}/5\mathbb{Z}\times \mathbb{Z}/10\mathbb{Z}$. One can see the extra $2$-torsion factor directly from Theorem \ref{thmain} as the polynomial $a(a_6x^6+a_4x^4+a_2x^2+a_0)$ is easily checked to have the quadratic factor $zx^2-1$. 
\end{remark}

\section{Search for a genus $2$ curve with a rational Weierstrass point and $\rk_5\Jac(C)(\mathbb{Q})_{\tors}=2$}

We ran a na\"ive search for solutions to \eqref{disceqn}, and for each triple of solutions $(t,u,z)\in\mathbb{Q}^3$ satisfying the hypotheses of Theorem \ref{thmain}, we computed the hyperelliptic curve of Theorem \ref{thmain} and checked if the curve $C$ had a rational Weierstrass point. For $t$ and $u$ rational numbers of the form $a/b$ with $|a|,|b|\leq 100$, we found 16630 appropriate solutions $(t,u,z)\in\mathbb{Q}^3$, which led to $274$ genus $2$ curves with a rational Weierstrass point and $\rk_5\Jac(C)(\mathbb{Q})=2$, lying in $85$ distinct isomorphism classes over $\mathbb{Q}$. For instance, the genus $2$ curve corresponding to the solution 
\begin{align*}
(t,u,z)=\left(\frac{2}{3},-\frac{1}{3},25\right), 
\end{align*}
may be given by the Weierstrass equation (after some simplification and changing to an odd model to make the rational Weierstrass point obvious):
\begin{align*}
C_0:y^2&=640x^5 + 3641x^4 + 8878x^3 + 11729x^2 + 8392x + 2576\\
&=(5x + 7)(128x^4 + 549x^3 + 1007x^2 + 936x + 368).
\end{align*}

The computer algebra system Magma \cite{Magma} verifies that the rational torsion subgroup of $\Jac(C_0)$ is exactly $\mathbb{Z}/5\mathbb{Z}\times \mathbb{Z}/10\mathbb{Z}$ (see Remark \ref{2rem}). Then the existence of the curve $C_0$ (or any of the other $84$ found curves) combined with Theorem \ref{chyp} finishes the proof of Theorem \ref{BLT}.

%\bibliographystyle{amsplain}
%\bibliography{5torsion}

\providecommand{\bysame}{\leavevmode\hbox to3em{\hrulefill}\thinspace}
\providecommand{\MR}{\relax\ifhmode\unskip\space\fi MR }
\providecommand{\MRhref}[2]{%
  \href{http://www.ams.org/mathscinet-getitem?mr=#1}{#2}
}
\providecommand{\href}[2]{#2}

\end{document}